\documentclass[12pt]{amsart}
\usepackage{amscd, amssymb,latexsym,amsmath, amscd, amsmath}
\usepackage[all]{xy}

\newcommand{\Q}{{\mathbb Q}}

\newcommand{\Z}{{\mathbb Z}}
\newcommand{\C}{{\mathbb C}}

\newcommand{\A}{{\mathbb A}}

\newtheorem{theorem}{Theorem}%[section] (If you want theorem numbered
\newtheorem{corollary}{Corollary}%       goes for lemmas, etc.)
\newtheorem{proposition}[theorem]{Proposition}
\theoremstyle{definition}
\newtheorem*{conjecture}{Conjecture}
\theoremstyle{remark}
\newtheorem*{remark}{Remark}
\newtheorem*{claim}{Claim}
\begin{document}

\title[Distinguished representations]
{Distinguished representations, base change, and reducibility for
 unitary groups}
\author{U.~K. Anandavardhanan}
\author{C.~S. Rajan}

\address{Tata Institute of Fundamental 
Research, Homi Bhabha Road, Bombay - 400 005, INDIA.}
\email{anand@math.tifr.res.in, rajan@math.tifr.res.in}


\date{}

\begin{abstract}

We show the equality of the local Asai $L$-functions 
defined via the Rankin-Selberg method and the Langlands-Shahidi method 
for a square integrable representation of $GL_n(E)$. As a consequence we
characterise reducibility of certain induced representations of 
$U(n,n)$, and the image of the base change map from $U(n)$ to $GL_n(E)$ in
terms of $GL_n(F)$-distinguishedness.    
\end{abstract}

\maketitle

\section{Introduction}
A representation $(\pi, V)$ of a group $G$ is said to be distinguished 
with respect to a character $\chi$ of a subgroup $H$,  
if there exists a linear form $l$ of $V$ satisfying
$l(\pi(h)v)=\chi(h)l(v)$ for all $v\in V$ and $h\in H$. When the
character $\chi$ is taken to be the trivial character, such
representations are also called as distinguished representations of
$G$ with respect to $H$. The concept of distinguished representations
can be carried over to a continuous context of representations of real
and  $p$-adic Lie groups, as well in a global automorphic context
(where the requirement of a non-zero linear form is replaced by the
non-vanishing of a period integral). 
The philosophy, due to Jacquet, is that 
 representations of a group $G$ distinguished with respect to 
a subgroup $H$ of fixed points of an involution on $G$ are often functorial
lifts from another group $G^\prime$. 

In this paper we consider 
$G={\rm Res}_{E/F}GL(n)$ and $H=GL(n)$ where
$E$ is a quadratic extension
 of a non-Archimedean local field $F$ of characteristic zero.  
In this case, the group $G^\prime$ 
is conjectured to be the quasi-split unitary group with respect
to $E/F$, 
\[G^\prime= U(n)=\{g \in GL_n(E) \mid gJ{^t\bar{g}}=J\},\]
where $J_{ij}=(-1)^{n-i}\delta_{i,n-j+1}$ and $\bar{g}$ is  the Galois
conjugate of $g$. There are two base change maps from $U(n)$ to
$GL(n)$ over $E$ called the stable and the unstable base change maps
(see Section \ref{section:bc}).  
We have the following conjecture due to Flicker and Rallis \cite{flicker2}:

\begin{conjecture} 
Let $\pi$ be an irreducible admissible representation of
$GL_n(E)$. If $n$ is odd (resp. even), then $\pi$ is $GL_n(F)$-distinguished 
if and only if it is
a stable (resp. unstable) base change from $U(n)$. 
\end{conjecture}

When $n=1$ the above conjecture is just Hilbert Theorem 90. The case $n=2$ is
established by Flicker \cite{flicker2}. 
The following theorem  proves the
 conjecture for a supercuspidal representation
when $n=3$. 
\begin{theorem} \label{bc}
A supercuspidal representation $\pi$ of $GL_3(E)$ is distinguished with respect
to $GL_3(F)$ if and only if it is a stable base change lift from $U(3)$.
\end{theorem} 

Let $G$ be a reductive $p$-adic group. 
An irreducible tempered representation of $G$
occurs as a component of an  induced
representation $I(\pi)$, parabolically induced  from a square integrable
representation $\pi$ of the Levi component $M$ of a parabolic
subgroup $P$ of
$G$. Thus the tempered spectrum of $G$ is determined from a knowledge
of the discrete series representations of the Levi components of
different parabolics and knowing the decomposition of induced
representations. The decomposition of $I(\pi)$ is governed by the
theory of $R$-groups. 

Let $G=U(n,n)$ be the quasi-split unitary group in $2n$
variables over
a $p$-adic field $F$, defined with respect to a quadratic extension
$E$ of $F$. Let $P$ be a parabolic subgroup of $G$ with a Levi
component $M$ isomorphic to $GL_{n_1}(E)\times  \cdots\times
GL_{n_t}(E)$ for some integers $n_i\geq 1$  satisfying $\sum_{i=1}^tn_i=n$.
Let $\pi_i, ~1\leq i\leq t$ be discrete series representations of
$GL_{n_i}(E)$. Let  $\pi=\pi_1\otimes\cdots\otimes\pi_t$ be
the discrete series representation of $M$. 
Let $\omega_{E/F}$ denote the quadratic character of $F^*$ associated to
the quadratic extension $E/F$.
The following theorem  gives
a description of the $R$-group $R(\pi)$ in terms of distinguishedness of the
representations $\pi_i$: 

\begin{theorem} \label{ind}
With the above notation, 
\[ R(\pi)\simeq (\Z/2\Z)^r,\]
where $r$ is the number of inequivalent representations $\pi_i$ 
 which are $\omega_{E/F}$-distinguished with respect to
$GL_{n_i}(F)$. 
\end{theorem}
\begin{corollary}
Let $P$ be  a maximal parabolic of $U(n,n)$ with Levi
component isomorphic to $GL_n(E)$, and $\pi$ be  a discrete series
representation of $GL_n(E)$.  Then $I(\pi)$ is 
reducible if and only if $\pi$ is $\omega_{E/F}$-distinguished
with respect to $GL_n(F)$.
\end{corollary}
 
A particular consequence of the corollary  is the 
following result about the 
Steinberg representation of $GL_n(E)$, which  is part of a more general 
conjecture, due to D. Prasad, about the Steinberg representation of 
$G(E)$ where $G$ is a 
reductive algebraic group over $F$ \cite{dprasad}.

\begin{theorem} \label{steinberg}
Let $\pi$ be the Steinberg representation of $GL_n(E)$.  Then $\pi$
is distinguished with respect to a character $\chi\circ \det$, for a character
$\chi$ of $F^*$, of $GL_n(F)$ if
and only if $n$ is odd and $\chi$ is the trivial character, or 
$n$ is even and $\chi=\omega_{E/F}$.
\end{theorem}

Our approach to the above theorems is via the theory of Asai
$L$-functions.  
The Asai $L$-function, 
also called the twisted tensor $L$-function,
 can be defined in three different ways: one via
the local Langlands correspondence and in terms of Langlands
parameters denoted by $L(s, As(\pi))$; via the theory of Rankin-Selberg
integrals \cite{flicker1,flicker3,kable} denoted by $L_1(s, As(\pi))$;
and the Langlands-Shahidi method (applied to a suitable unitary 
group) \cite{goldberg,shahidi2} denoted by $L_2(s, As(\pi))$. 
 It is of course expected that all the above three 
$L$-functions match. 

The main point is that the analytical properties of the different
definitions of Asai $L$-function give different insights about the
representation: the Asai $L$-function defined via the 
Rankin-Selberg method
can be related to distinguishedness with respect to $GL_n(F)$, whereas
the Asai $L$-function defined via the Langlands-Shahidi method is related
to the base change theory from $U(n)$, and to reducibility questions for
$U(n,n)$. 
Thus the following theorem, proved using global methods, 
is a key ingredient towards a proof
of the above theorems: 
\begin{theorem}\label{matching}
Let $\pi$ be a square integrable representation of $GL_n(E)$. Then
$L_1(s,As(\pi))=L_2(s,As(\pi))$.
\end{theorem}      

\section{Asai $L$-functions}
\subsection{Langlands parameters}  Let $F$ be a
non-archimedean local field and let $E$ be  a quadratic
extension of $F$. The  Weil-Deligne group $W^\prime_E$ of $E$ is of
index two in the Weil-Deligne group $W^\prime_F$ of $F$.  Choose 
$\sigma \in W^\prime_F\backslash W^\prime_E$ of order $2$. 
Given a  a continuous, $\Phi$-semisimple
representation $\rho$ of $W^\prime_E$ of dimension $n$, the
representation $As(\rho): W^\prime_F\to GL_{n^2}(\C)$ given by tensor
induction of $\rho$ is defined as, 
\[ As(\rho)(x)=\begin{cases} \rho(x)\otimes \rho(\sigma^{-1} x \sigma)
&\text{if $x\in W^\prime_E$}\\
\rho(\sigma x)\otimes \rho(x\sigma) &\text{if  $x\not \in W^\prime_E$}.
\end{cases}
\]
Let $\pi$ be an irreducible, admissible representation of $GL_n(E)$
with Langlands parameter $\rho_{\pi}$. The Asai $L$-function $L(s,
As(\pi))$ is defined to be the $L$-function $L(s, As(\rho_{\pi}))$. 

\subsection{Rankin-Selberg method}
\subsubsection{Local theory}
We recall the Rankin-Selberg theory of the Asai $L$-function 
\cite{flicker1}, \cite{flicker3}, and \cite{kable}. Let $F$ be a
non-archimedean local field and let $E$ be either a quadratic
extension of $F$ or $F\oplus F$.  
Let $\pi$ be an irreducible admissible generic representation of $GL_n(E)$.
We take an additive character $\psi$ of $E$ which restricts trivially
to $F$. There exists an additive character $\psi_0$ of $F$ such that
$\psi(x)=\psi_0(\Delta(x-\bar{x}))$ where $\Delta$ is a trace zero element 
of $E^*$.
Let ${\mathcal W}(\pi,\psi)$ denote the  Whittaker model of $\pi$ 
 with respect to $\psi$.
Let $N_n(F)$ be the unipotent radical of the Borel subgroup of $GL_n(F)$.
Consider the integral (see \cite{flicker1})
$$\Psi(s,W,\Phi)=\int_{N_n(F)\backslash GL_n(F)}W(g)\Phi((0,0,\ldots,1)g)
|\det g|_F^sdg,$$
where $\Phi\in {\mathcal S}(F^n)$, the space of locally constant compactly 
supported functions on $F^n$, and $dg$ is a $GL_n(F)$-invariant measure on 
$N_n(F)\backslash GL_n(F)$.

In \cite{flicker3}, Flicker proves that the above integral converges
absolutely in some right half plane to a rational function in $X=q^{-s}$,
where $q=q_F$ is the cardinality of the residue field of $F$. 
The space spanned by
$\Psi(s,W,\Phi)$ (as $W$ and $\Phi$ vary) is a fractional ideal in 
${\Bbb C}[X,X^{-1}]$ containing the constant function $1$. We can choose
a unique generator of this ideal of the form $P_1(X)^{-1}$, $P_1(X)
\in {\Bbb C}[X]$
such that $P_1(0)=1$. Define the  Asai $L$-function $L_1(s,As(\pi))$
as 
\[ L_1(s,As(\pi))=P_1(q^{-s})^{-1}.\] 
This does not depend on the choice of the additive character
$\psi$.  Moreover $\Psi(s,W,\Phi)$ satisfies the functional equation
$$\Psi(1-s,\widetilde{W},\hat{\Phi})=\gamma_1(s,As(\pi),\psi)\Psi(s,W,\Phi)$$
where $\widetilde{W}(g)=W(w{^tg^{-1}})$, 
$w$ is the longest  element of the Weyl group,
and $\hat{\Phi}$ is the Fourier transform of $\Phi$ with respect to $\psi_0$.
The epsilon factor, 
$$\epsilon_1(s,As(\pi),\psi)=\gamma_1(s,As(\pi),\psi)\frac{L_1(s,As(\pi))}
{L_1(1-s,As(\pi^\vee))}$$  is a monomial
in $q_F^{-s}$.

If  $E=F\oplus F$ write $\pi=\pi_1\times \pi_2$ considered as a
representation of $G_n(F)\times GL_n(F)$. Then
\begin{equation}\label{rs-split}
 L_1(s, As(\pi))= L(s, \pi_1\times \pi_2),
\end{equation}
where the right hand side is  the Rankin-Selberg $L$-factor of
$\pi_1\times \pi_2$. 

We  have the following proposition \cite[Proposition in
  Section 3]{flicker1}:
\begin{proposition}\label{unram-rs}
Suppose $E/F$ is an unramified quadratic extension.
Let $\pi=Ps(\mu_1,\ldots,\mu_n)$ be an unramified unitary representation
induced from the character $(t_1,\ldots,t_n)\longrightarrow \prod \mu_i(t_i)$
of the diagonal torus in $GL_n(E)$.
Let $W_{\pi}^0$ be  the spherical Whittaker
function, and 
$\Phi_{F}^0$ be  the characteristic function of ${\mathcal O}_F^n$.
Then 
\[\Psi(s,W_{\pi}^0,\Phi_{F}^0)
=\prod_{j=1}^n(1-\mu_j(\varpi_F)q_F^{-s})^{-1}\cdot
\prod_{i<j}(1-\mu_i(\varpi_F)\mu_j(\varpi_F)q_F^{-2s})^{-1}\]
where $\varpi_F$ is a uniformizing parameter of $F$.
\end{proposition}

The following proposition is proved in \cite[Theorem 4]{kable}:
\begin{proposition} \label{js-rs}
Let $\pi$ be a square integrable representation of $GL_n(E)$. Then
$L_1(s,As(\pi))$ is regular in the region ${\rm Re}(s)>0$. 
\end{proposition}
We remark that for the proof of Theorem \ref{matching} all that we
require is that $L_1(s,As(\pi))$ be regular in the region ${\rm
  Re}(s)\geq 1/2$.

\subsubsection{Global theory}
Now let $L/K$ be a quadratic extension of number fields. We assume that
the archimedean places of $K$ split in $L$. Let $\psi_0$
be a non-trivial character of ${\Bbb A}_K/K$, and let $\psi=\psi_0
(\Delta(x-\bar{x}))$.  
 For a global
field $K$, let $\Sigma_K$ denote the set of places of $K$. Let $\Pi=
\bigotimes_{w\in \Sigma_L} \Pi_w$ be a 
 representation of $GL_n({\Bbb A}_L)$.
Let $T$ be a finite set of places of $K$ containing the following
places:
\begin{itemize}
\item the archimedean places of $K$. 
\item the ramified places of  the extension $L/K$.
\item the places $v$ of $K$ dividing a  place  $w$ of $L$, where
  either $\psi_{0,v}, ~\psi_{L_w}$ or $\Pi_w$ is ramified.  
\end{itemize}

Define, 
\begin{equation}
L^\prime_{1,v}(s,As(\Pi))=\begin{cases} L_1(s, As(\Pi_w)) & \text{ $w|v$,
 $v\in T$ and 
 $v$ inert},\\
\Psi_v(s, W_{\Pi_w}^0, \Phi_{F_v}^0) & \text{$v$ inert, $v\not\in T$},\\
L(s, \Pi_{w_1}\times \Pi_{w_2}) &\text{$v$ splits, $v=w_1w_2$}.
\end{cases}
\end{equation}
\begin{remark} Let $v$ be a place of $K$ not in $T$, inert in $L$ and
 $w$ the place of $L$ dividing $v$.  
 It is not known that  $L_1(s,As(\Pi_w))=\Psi(s,
W_{\Pi_w}^0, \Phi_{K_v}^0)$. In the notation of
Proposition \ref{unram-rs}, the right
hand side is the $L$-factor associated by Langlands functoriality. 
\end{remark}

Following Kable \cite{kable}, we define the Rankin-Selberg Asai
$L$-function $L_1(s,As(\Pi))$ as, 
\[L_1(s,As(\Pi), T)=\prod_{v\in \Sigma_K} L^\prime_{1,v}(s,As(\Pi)).
\]
We have the functional equation:
\begin{proposition}[Theorem 5, \cite{kable}] \label{fe-rs}
Let $\Pi$ be a cuspidal automorphic representation of
$GL_n(\A_L)$. Then $L_1(s,As(\Pi), T)$ admits a meromorphic
continuation to the entire plane and satisfies the functional equation
\[ L_1(s,As(\Pi), T)= \epsilon_1(s,As(\Pi),T)L_1(1-s,As(\Pi^\vee), T)\]
where  the function $\epsilon_1(s,As(\Pi),T)$ is entire and non-vanishing. 
\end{proposition}

\subsection{Langlands-Shahidi method}
\subsubsection{Local theory}
We now recall the  Langlands-Shahidi 
approach to the Asai $L$-function \cite{goldberg, shahidi2}.
Let 
$ G=U(n,n)$ be the quasi-split unitary group in 
$2n$ variables with respect to $E/F$. The group $M=R_{E/F}GL_n$
can be embedded as a Levi component of a maximal parabolic subgroup
$ P$ of $ G$ with unipotent radical $N$.  Let $r$ be the
adjoint representation of the $L$-group of $M$ on the Lie algebra of 
the $L$-group of $N$. Fix an additive character $\psi_0$ of $F$.
The Langlands-Shahidi gamma factor $\gamma_2(s,\pi,r,\psi_0)$ 
defined in \cite{shahidi2}, is a  rational function of  $q^{-s}$.
Let $P_2(X)$ be the unique polynomial satisfying $P_2(0)=1$
such that $P_2(q^{-s})$ is the numerator of $\gamma_2(s,\pi,r,\psi_0)$.
For a tempered $\pi$, the Langlands-Shahidi Asai $L$-function is
defined as
\[L_2(s,As(\pi))=1/P_2(q^{-s}).\]
The $L$-function is independent of the additive character. 
The quantity  
\[ \epsilon_2(s,As(\pi),\psi_0)=\gamma_2(s,\pi,r,\psi_0)\frac{L_2(s,As(\pi))}
{L_2(1-s,As(\pi^\vee))}\] is the
Langlands-Shahidi epsilon factor, and is a monomial in $q^{-s}$. 

The analytical properties of $L_2(s,As(\pi))$ are 
proved in \cite[Theorem 3.5, Proposition 7.2]{shahidi2}:
\begin{proposition} \label{local-ls}
Let $\pi$ be an irreducible admissible  representation of
$GL_n(E)$. Then the following holds:
\begin{enumerate}
\item If $E$ is an unramified extension of $F$, and 
$\pi=Ps(\mu_1,\ldots,\mu_n)$ be a  unitary unramified representation
  of $GL_n(E)$, as in the hypothesis of Proposition \ref{unram-rs}. 
Then 
\[L_2(s,As(\pi))=\prod_{j=1}^n(1-\mu_j(\varpi_F)q_F^{-s})^{-1}\cdot
\prod_{i<j}(1-\mu_i(\varpi_F)\mu_j(\varpi_F)q_F^{-2s})^{-1}.\]

\item  Let $\pi$ be a tempered  representation of $GL_n(E)$. 
Then
$L_2(s,As(\pi))$ is regular in the region ${\rm Re}(s)>0$. 
\end{enumerate}
\end{proposition}

\subsubsection{Global theory}
Let $L/K$ be a quadratic extension of number fields, and let $\Pi=
\bigotimes_w \Pi_w$ be a
 representation of $GL_n({\Bbb A}_L)$.   Define for a place $v$ of
 $K$,  
\begin{equation}
L_{2,v}(s,As(\Pi))=\begin{cases} L_2(s, As(\Pi_w)) & \text{ $w|v$,
 $v$ inert},\\
L(s, \Pi_{w_1}\times \Pi_{w_2}) &\text{$v$ splits, $v=w_1w_2$}.
\end{cases}
\end{equation}
Define the global $L$-function 
$$L_2(s,As(\Pi))=\prod_{v\in \Sigma_K} L_{2,v}(s,As(\Pi)).$$ 
Then we have the functional equation \cite{shahidi2}:
\begin{proposition} \label{fe-ls}
Let $\Pi$ be a cuspidal automorphic representation of
$GL_n(\A_L)$. Then $L_2(s,As(\Pi))$ admits a meromorphic
continuation to the entire plane and satisfies the functional equation
\[ L_2(s,As(\Pi))=\epsilon_2(s,As(\Pi))L_2(1-s,As(\Pi^\vee))\]
where  the function $\epsilon_2(s,As(\Pi))$ is entire and non-vanishing. 
\end{proposition}

\section{Proof of Theorem \ref{matching}}
The proof of Theorem \ref{matching} is via global methods. The
following proposition 
embedding a square integrable representation $\pi$ as the local
component of a cuspidal automorphic representation is well known
\cite[Lemma 5]{kable}, \cite[Lemma 6.5 of Chapter 1]{arthur}:
\begin{proposition} \label{embed}
Let $E/F$ be a quadratic extension of non-archimedean local fields of
characteristic zero and residue characteristic $p$. Let $\pi$ be a square
integrable representation of $GL_n(E)$. Then the following holds: 
\begin{enumerate}
\item There exists a number field $K$, a quadratic extension $L$ of
  $K$ and a place $v_0$ of $K$ inert in $L$, such that $K_{v_0}\simeq F$ and
  $L_{w_0}\simeq E$, where $w_0$ is the unique place of $L$ dividing
  $v_0$. Further, $v_0$ is the unique place of $K$
  lying over the rational prime $p$, and the real places of $K$ are
  split in $L$.  
\item There exists a cuspidal automorphic representation $\Pi$ of 
$GL_n(\A_L)$ such that $\Pi_{w_0}\simeq \pi$. 
\end{enumerate}
\end{proposition}
Let $\Pi$ be a cuspidal representation of $GL_n(\A_L)$ satisfying the
properties of the above proposition. Choose
 a finite set $T$ of places
of $K$ as in  Proposition \ref{fe-rs}.  Consider the ratio 
\[F(s, \Pi)=\frac{L_2(s,As(\Pi))}{L_1(s,As(\Pi),T)}.\]
If $v=w_1w_2$ is a place of $K$ which splits into two places $w_1$ and
$w_2$ of $L$, then 
\[ L^\prime_{1,v}(s,As(\Pi))=L_{2,v}(s,As(\Pi))=
L(s, \Pi_{w_1}\times \Pi_{w_2}) .\]
By
Propositions \ref{unram-rs} and \ref{local-ls}, 
if $v$ is a place of $K$ which is inert and not in $T$, then 
\[ L^\prime_{1,v}(s,As(\Pi))=L_{2,v}(s,As(\Pi)).\]
Hence 
\[F(s, \Pi)=\prod_{v\in T}
\frac{L_{2,v}(s,As(\Pi))}{L^\prime_{1,v}(s,As(\Pi))}. \]
Write 
\[F(s, \Pi)=G(s, \Pi)Q(s,\Pi)P_0(s,\Pi),\]
where 
\begin{itemize}
\item The function $G(s,\Pi)$ is the ratio of the $L$-factors at the
  archimedean places; it  is a ratio of products of Gamma functions of the form
  $\Gamma(as+b)$ for some suitable constants $a, ~b$. 
\item The function
\[Q(s,\Pi)=\frac{\prod_{i=1}^n (1-{\alpha_i}q_{v_i}^{-s})}
{\prod_{j=1}^m (1-{\beta_j}q_{v_j}^{-s})}, \quad v_i, ~v_j\in
  T':=T\backslash \{v_0\}\]
is a ratio of the $L$-factors at the
  finite set of places
  of $T$ not equal to $v_0$; it is a ratio of 
products of distinct functions of the form $(1-{\beta}
  q_v^{-s}), ~\beta\neq 0$, 
 where $v\in T':=T\backslash \{v_0\}$, and $q_v$ is the number
  of elements of the residue field. By our assumption on $K$, $(p,
  q_v)=1$. 
\item  The function  
\[P_0(s,\Pi)=\frac{L_2(s, As(\pi))}{L_1(s,
  As(\pi))}\]
 is a  ratio of products of functions of the form $(1-{\alpha}
  q_{v_0}^{-s})$. 
  
By Propositions \ref{js-rs} and \ref{local-ls}, the functions
  $P_0(s,\Pi)$ and $P_0(s,\Pi^\vee)$ are 
  regular and non-vanishing in the region ${\rm Re}(s)\geq 1/2$. 
\end{itemize}

We claim the following: 

\noindent{\em Claim.} Let $\gamma_0$ be a pole (resp. zero) 
of $P_0(s, \Pi)$. The function $F(s,\Pi)$ has a pole (resp. zero) at all
but finitely many elements of the 
form $\gamma_0+2\pi ik/{\rm log}~q_{v_0}$, $k\in \Z$.

\begin{proof}[Proof of Claim] Suppose   that the function
$F(s,\Pi)$ is regular at  points of the form
$\gamma_0+2\pi il/{\rm log}~q_{v_0}$ for  integers $l\in C$, where $C$
is an infinite subset of the integers. 
 Since $G(s)$ can contribute only finitely many zeros on any line
with real part constant, these poles have to be cancelled by zeros of
$Q(s,\Pi)$. Since $T$ is finite, and the local $L$-factors are
polynomial functions in $q_v^{-s}$, there is a $v\in T', ~\gamma\in
\C$  and a function $f:C\to \Z$ such that,
\[ \gamma_0+2\pi il/{\rm log}~q_{v_0}=\gamma+2\pi if(l)/{\rm
  log}~q_{v}\]
for infinitely many $l \in C$.
Taking the difference of any two elements, we get 
${\rm log}~q_{v_0}/{\rm log}~q_{v}\in \Q$. 
This is not possible as
$q_{v_0}$ and $q_v$ are coprime integers. Hence, all but finitely
many poles of the form $\gamma_0+2\pi ik/{\rm log}~q_{v_0}, ~k\in \Z$
are poles of
$F(s,\Pi)$. 
\end{proof}

Since $P_0(s,\Pi)$ is
  regular in the region ${\rm Re}(s)\geq 1/2$, we obtain
${\rm Re}(\gamma_0)<1/2$. 
From the global functional equations given by Propositions \ref{fe-rs} and
\ref{fe-ls}, $F(s, \Pi)$ satisfies a functional equation, 
\[ F(s,\Pi)=\eta(s,\Pi)F(1-s,\Pi^\vee),\]
where $\eta(s,\Pi)$ is an entire non-vanishing function. 
Hence $F(s, \Pi^\vee)$ has
infinitely many poles of the form $1-\gamma_0+2\pi ik/{\rm
  log}~q_{v_0}$ with $k \in \Z$. Since  $P_0(s,\Pi^\vee)$ is
  regular in the region ${\rm Re}(s)\geq 1/2$, these poles have to be
  poles of $G(s,\Pi^\vee)Q(s,\Pi^{\vee})$. Arguing as in proof of the above
  claim, we obtain a contradiction. 
Arguing similarly with the zeros instead of poles, we obtain that
$P_0(s,\Pi)$
is an entire non-vanishing function 
and hence it is a constant. Since the
  $L$-factors are normalised, we obtain a proof of Theorem \ref{matching}. 

\begin{remark} The method of proof of Theorem \ref{matching} is a
  general method allowing us to establish an equality for two
  possibly different definitions of $L$-factors at `bad' places. This
  requires a global functional equation, equality of the $L$-factors
  at all good places, and regularity in the region ${\rm Re}(s)\geq 1/2$
for the `bad' $L$-factors. The method is illustrated in \cite{rajan} in
  the context of functoriality, but allowing the use of cyclic base change. It
  is used by Kable in \cite{kable} to prove, for a square integrable
  representation, that the Rankin-Selberg $L$-factor 
$L(s,\pi\times \bar{\pi})$ factorises as a product of $L_1(s, As(\pi))$
  times $L_1(s, As(\pi \otimes \widetilde{\omega}))$, where 
$\widetilde{\omega}$ is an extension of $\omega_{E/F}$, the
  quadratic character corresponding to the extension $E/F$. A proof of
  strong multiplicity one in the Selberg class using similar arguments
  is given in \cite{murty}. 
\end{remark}

\begin{remark} It has been shown by Henniart \cite{henniart} using
  similar global methods, that for
  any irreducible, admissible representation $\pi$ of $GL_n(E)$, the
  equality  
$L(s, As(\pi))=L_2(s, As(\pi))$.  Henniart's
  proof uses cyclic base change and the inductivity of
  $\gamma$-factors  to go from square integrable to all irreducible,
  admissible representations.  Since we do not know inductivity of the
  Rankin-Selberg $\gamma$-factors $\gamma_1(s, As(\pi), \psi)$, we
  cannot derive a similar statement for the Rankin-Selberg $L$-factors.  
\end{remark}

\begin{remark} Using cyclic base change as in \cite{rajan} or
  \cite{henniart}, it is possible to show that the $\epsilon$-factors
$\epsilon_1(s, As(\pi), \psi)$ and $\epsilon_2(s, As(\pi), \psi_0)$
  are equal up to a root of unity, when $\pi$ is square integrable. 
\end{remark}

\section{Applications}

\subsection{Analytic characterisation of distinguished
  representations} 

The proofs of Theorems \ref{bc} and \ref{ind} use the 
following proposition relating  the concept of distinguishedness
with the analytical properties of the (Rankin-Selberg) Asai
$L$-function \cite[Corollary 1.5]{anand}:
\begin{proposition}\label{dist}
Let $\pi$ be a square integrable representation of $GL_n(E)$. Then
$\pi$ is distinguished with respect to $GL_n(F)$ if and only if 
$L_1(s,As(\pi))$ has a pole at  $s=0$. 
\end{proposition}

\subsection{Base change for $U(3)$} \label{section:bc}

 Let $W_{E/F}$ be the relative Weil group of $E/F$ defined as
the  semidirect product of $E^*\rtimes {\rm Gal}(E/F)$ for the natural
action of  ${\rm Gal}(E/F)$ on $E^*$. The Langlands dual group of
$U(n)$ is given by 
$^LU(n)=GL_n(\C)\rtimes W_{E/F}$, where $W_{E/F}$ acts via the
projection to  ${\rm Gal}(E/F)$, and the non-trivial element $\sigma
\in {\rm Gal}(E/F)$ acts by $\sigma(g)=J~^tg^{-1}J^{-1}$ on
$GL_n(\C)$. The Langlands dual group of $R_{E/F}(GL_n)$ is given by 
\[^LR_{E/F}(U(n)) =[GL_n(\C)\times GL_n(\C)]\rtimes W_{E/F}.\]
Here again the action of  $W_{E/F}$ is via the
projection to  ${\rm Gal}(E/F)$, and $\sigma$ acts by 
$(g,h)\mapsto (J~^th^{-1}J^{-1},J~^tg^{-1}J^{-1})$. 

There are two natural mappings from the $L$-group of $U(n)$ to the $L$-group of
$GL_n(E)$, called the stable and the unstable base change maps.
At the $L$-group
level, the stable base change map, which corresponds to the restriction
of parameters from the Weil group $W_F$ of $F$ to the Weil group $W_E$
of $E$, is given by the diagonal
embedding $\psi:{^LU(n)}\to {^LR_{E/F}(U(n))}$. The unstable base change
map is defined by first choosing a character $\widetilde{\omega}$ of
$E^*$ extending the quadratic character $\omega_{E/F}$ of $F^*$ associated to
the quadratic extension $E/F$. At the level of $L$-groups, the
unstable base change corresponds to the homomorphism 
$\psi^\prime: {^LU(n)}\to {^LR_{E/F}(U(n))}$ 
given by $\psi^\prime(g\times w)
=(\widetilde{\omega}(w)g, \widetilde{\omega}(w)^{-1}g)\times
w$ for $w\in
E^*, ~g\in GL_n(\C)$ and $\psi^\prime(1,\sigma)=(1, -1)\times
\sigma$. The base change lift for $n=3$ has
been established by Rogawski \cite{rogawski}.

\begin{proof}[Proof of Theorem \ref{bc}] 
By \cite[Corollary 4.6]{goldberg}, a supercuspidal representation 
$\pi$  of $GL_3(E)$ is a stable
  base change lift from $U(3)$ if and only if 
 the Langlands-Shahidi Asai $L$-function $L_2(s,As(\pi))$
has a pole at $s=0$. By Theorem \ref{matching}, this amounts to 
saying that the Rankin-Selberg  Asai $L$-function $L_1(s,As(\pi))$
has a pole at $s=0$. Now the theorem follows by appealing to
Proposition \ref{dist}.
\end{proof}

\begin{remark}
If $\pi$ is a square integrable representation such that 
$\pi^\vee \cong \bar{\pi}$,
and the central character of $\pi$ has trivial restriction to $F^*$,
 then Kable \cite{kable} has proved that
$\pi$ is distinguished or distinguished with respect to $\omega_{E/F}$,
the quadratic character associated to the extension $E/F$ (see 
\cite{hakim2,dprasad} for earlier results in this direction). 
The given conditions on $\pi$ are expected to be necessary  
for $\pi$ to be
in the image of the base change map from $U(n)$. Thus
Kable's result can be thought of as a weaker version of the conjecture stated
in the introduction.     
On the other hand, it is expected that $U(n)$-distinguished representations
of $GL_n(E)$ are base change lifts from $GL_n(F)$. This has been proved
in several cases \cite{hakim1,dprasad}. 
\end{remark}

\subsection{Reducibility for $U(n,n)$}

We now prove  Theorem \ref{ind}. 
 In \cite{goldberg,goldberg1}, Goldberg proves that
for a discrete series representation $\pi$ with
$\pi^\vee \cong \bar{\pi}$, $I(\pi)$ is irreducible if and only
if $L_2(s,As(\pi))$ has a pole at $s=0$ (see also \cite{ichino}). 
By \cite[Theorem 3.4]{goldberg1}, $R(\pi)\simeq (\Z/2\Z)^r$, where $r$
is the number of inequivalent representations $\pi_i$ satisfying
$\pi_i^\vee\simeq \bar{\pi}_i$ and the Plancherel measure
$\mu(s,\pi_i)$ has no zero at $s=0$. By \cite[Corollary 3.6]{shahidi2}, 
the latter condition amounts to knowing that the Asai
$L$-functions $L_2(s, As(\pi_i))$ are regular at $s=0$. 

Theorem \ref{ind} follows from the following claim:
\begin{claim} An irreducible, square integrable representation $\pi$
  of $GL_n(E)$ is $\omega_{E/F}$ distinguished if and only if 
$\pi^\vee\simeq \bar{\pi}$ and $L_2(s, As(\pi))$ is regular at $s=0$.
\end{claim}

\noindent
{\em Proof of Claim.}
By \cite[Corollary 5.7]{goldberg}, 
\[ L(s,\pi\times \bar{\pi})=L_2(s, As(\pi))L_2(s,
As(\pi\otimes\widetilde{\omega})),\]
where $\widetilde{\omega}$ is a character of $E^*$ which
restricts to ${\omega_{E/F}}$ on $F^*$. Now $L(s,\pi\times \bar{\pi})$
has a pole at $s=0$ if and only if $\pi^\vee\simeq \bar{\pi}$. Hence
$\pi^\vee\simeq \bar{\pi}$ and $L_2(s, As(\pi))$ is regular at $s=0$
is equivalent to saying that $L_2(s,
As(\pi\otimes\widetilde{\omega}))$ has a pole at $s=0$. 
By Theorem \ref{matching} this is the same as saying that  $L_1(s,
As(\pi\otimes\widetilde{\omega}))$ has a pole at $s=0$. By
Proposition \ref{dist}, the latter condition is equivalent to
saying that $\pi$ is $\omega_{E/F}$ distinguished.
This proves Theorem \ref{ind}. 

\begin{remark}
The $R$-group in this context is also computed in terms of the Langlands
parameters by D. Prasad \cite[Proposition 2.1]{dprasad1}. 
According to this computation, $R(\pi)$ is a product of $r$ copies
of $\Z/2\Z$'s, where $r$ is the number of $\pi_i$'s such that 
$\pi_i^\vee \cong \bar{\pi_i}$, and $c(\sigma_i)=-1$, where $\sigma_i$
is the Langlands parameter of $\pi_i$. Here $c(\sigma_i)\in
\{\pm 1\}$
denotes the constant introduced by Rogawski \cite[Lemma 15.1.1]{rogawski}.
Also $c(\sigma_i)=(-1)^{n_i-1}$ if and only if $\sigma_i$ can be extended to
a parameter  for                                     
$U(n_i)$. Together with Theorem \ref{ind}, this gives an evidence
for the conjecture stated in the introduction.
\end{remark}

\subsection{Distinguishedness of Steinberg representation of $GL(n)$}  

Now let $G=GL(n)$. For a representation $\pi$ of $GL_n(E)$, let $I(\pi)$
be the parabolically induced representation of $U(n,n)$. If $\pi$ is a discrete
series representation such that $\pi^\vee
\ncong \bar{\pi}$, then $I(\pi)$ is known to be irreducible \cite{goldberg}.
Suppose $\pi^\vee
\cong \bar{\pi}$. Let $a$ and $b$ be integers such that $ab=n$, such that
$\pi$ is the unique square integrable constituent of the representation
induced from $\pi_1\otimes\ldots\otimes\pi_b$ where $\pi_i=\pi_0\otimes
|~|_E^{b+1-2i/2}$, and $\pi_0$ a supercuspidal 
representation of $GL_a(E)$. Then
$\pi_0^\vee \cong \bar{\pi_0}$. We have the following result
of Goldberg \cite[Section 7]{goldberg}:
\begin{proposition}
The representation $I(\pi)$ of $U(n,n)$ is irreducible if and only if
$L_2(s,As(\pi_0))$ (resp. $L_2(s,As(\pi_0\otimes \widetilde{\omega}))$ has 
a pole 
at $s=0$ if $b$ is odd (resp. even). Here $\widetilde{\omega}$ is a 
character of $E^*$
that restricts to $\omega_{E/F}$.
\end{proposition}
Now if $\pi$ is the Steinberg representation of $GL_n(E)$,
then $a=1$, $b=n$, and $\pi_0$ is the trivial character. Thus $I(\pi)$
is irreducible when $n$ is odd, and reducible when $n$ is even.
By the corollary to Theorem \ref{ind}, $\pi$ is $\omega_{E/F}$-distinguished
when $n$ is even, and $\pi$ is not $\omega_{E/F}$-distinguished when $n$ 
is odd. 

Since $\pi^\vee \cong \bar{\pi}$ and $\omega_\pi=1$, we know that $\pi$ is
either distinguished or $\omega_{E/F}$-distinguished, but not both 
(see \cite[Theorem 7]{kable} and \cite[Corollary 1.6]{anand}). Therefore
it follows that when $n$ is odd (resp. even), $\pi$ is distinguished 
(resp. $\omega_{E/F}$-distinguished), and that $\pi$ is not distinguished
with respect to any other character.

\section*{Acknowledgements} 
The first named author would like to thank Dipendra
Prasad for many helpful suggestions. We thank F. Shahidi and
G. Henniart for bringing to our notice the preprint \cite{henniart}.

\end{document}